\documentclass[numreferences]{kluwer}
\usepackage{graphicx}
\usepackage{setspace}
\usepackage{enumitem}
\usepackage{amsthm}

\newtheorem{Thm}{Theorem}[section]

\newtheorem{Alg}[Thm]{Algorithm}

\newtheorem{Rem}{Remark}

\begin{document}
	\begin{opening}
		\title{Trigonometric Interpolation Based Approach for Second Order ODE with Mixed Boundary Conditions}
		\author{Xiaorong \surname{Zou}\email{xiaorong.zou@bofa.com}}
		\institute{Global Market Risk Analytic, Bank of America}
		\runningauthor{X. Zou}
		\runningtitle{TIBA for Second Order ODE}
		\date{May 6, 2025}
		\classification{MSC2000}{Primary 65T40, Secondary 65T50}
		
		\begin{abstract}
			In this paper, we propose a  trigonometric-interpolation based approach (TIBA) to approximate solutions of mixed boundary value problems of
			second order ODEs.  TIBA leverages analytic attractiveness of a trigonometric polynomial to reformulate the dynamics of $y,y',y''$  implied by ODE and boundary conditions. TIBA is particularly attractive for a linear ODE where solution can be obtained directly by solving a linear system.  The framework can be used to solve integro-differential equations. Numerical tests  have been conducted to assess TIBA's performance regarding convergence, existence and uniqueness of solution under various boundary conditions with expected results. 			
		\end{abstract}
	\end{opening}

\section{Introduction}\label{sec:intro}

A new trigonometric interpolation algorithm was recently introduced in \cite{zou_tri}.  It leverages Fast Fourier Transformation (FFT) to achieve optimal computational efficiency and converges at speed aligned with smoothness of underlying function. In addition, it can be used to approximate nonperiodic functions defined on bounded intervals.  Considering the analytic attractiveness of trigonometric polynomial,  especially in handling differential and integral operations, the proposed trigonometric estimation of a general function is expected to be used in a wide spectrum. 
In this paper, we continue on applications of the trigonometric interpolation algorithm to solve the following nonlinear ODE system:
\begin{equation}\label{eq:nonlinear_ode_order2}  
	y''(x) = f(x,y,y'), \quad x\in [s,e]
\end{equation}
\begin{eqnarray}
	d_{11}y(s) + d_{12}y'(s) +d_{13}y(e) + d_{14}y'(e) &=& \alpha , \label{eq:nonlinear_ode_order2:diri}\\ 
	d_{21}y(s) + d_{22}y'(s) +d_{23}y(e) + d_{24}y'(e) &=& \beta, \label{eq:nonlinear_ode_order2:neum} 
\end{eqnarray}
where $f(x,v,u)$
is continuously differential on the range $[s, e]\times R^2$, the rank of matrix $D:=(d_{ij})_{1\le i\le 2, 1\le j\le 4}$ is $2$, and $\alpha, \beta$ are two real numbers. 

There are quite rich researches on algorithms of numerical solutions of boundary problems of second order ODE 
\cite{ode_naga}-\cite{ode_bumozh13}.
In \cite{zou_tri_III}, a trigonometric interpolation based optimization algorithm, labeled as $TIBO$, is proposed for the solutions of ODE  (\ref{eq:nonlinear_ode_order2}-\ref{eq:nonlinear_ode_order2:neum}). 
TIBO leverages analytic attractiveness of a trigonometric polynomial to discretize ODE (\ref{eq:nonlinear_ode_order2}) in a global manner (see Eq. (\ref{eq:bdp_orde2_d_nonlinear})). 
The boundary conditions can be captured by two parameters in a derived approximation of $y$, similar to how initial condition is captured in Adomian decomposition method \cite{ode_Adomian}.  In addition to its flexibility to address general non-linear ODEs with mixed boundary conditions, it can achieve high accuracy when optimization process converges, and address the issue of multiple solutions by integrating certain requirements in optimization to identify a desired solution. 

As a disadvantage of TIBO, for certain types of boundary conditions, its performance can be sensitive to initial guess of a solution at grid points used by optimization process.  For Neumann condition with give $y(s),y'(s)$,  we can use the standard Runge-Kutta (RK) scheme \cite{ryts} to generate initial values and thus $TIBO$ always generates accurate results.  For other types of condition, one can use some existing algorithm, such as shooting method \cite{HBKeller}, to approximate $y(s), y'(s)$, and then apply RK to generate initial values.  For Dirichlet condition with given $y(s)$,  shooting-RK combination can generate proper initial values and TIBO also generates descent results as in Neumann.  For a boundary condition where neither $y(s)$ nor $y'(s)$ is known, the algorithm becomes more sensitive to the initial guess of $y(s), y'(s)$ if shooting-RK combination is used to generate initial values and divergent scenarios are observed with significant chance as shown in \cite{zou_tri_III}.

In this paper,  we propose a similar trigonometric interpolation based algorithm, named TIBA hereafter.  
In a nutshell,   TIBA  discretizes globally ODE system (\ref{eq:nonlinear_ode_order2}-\ref{eq:nonlinear_ode_order2:neum}) into a non-linear system whose solutions can be solved by certain existing schemes like Newton method. TIBA is particularly attractive to solve a linear ODE by converting it to a linear algebraic system whose solution can be solved directly without requirement on initial guess about a solution as in TIBO.  With a linear system, TIBA outperforms TIBO in efficiency and addresses properly the issues of existence and uniqueness of solution based on derived linear algebraic system.  In addition,  TIBA can be extended to solve linear Integro-differential equations as studied in \cite{zou_tri_fredholm} and \cite{zou_tri_volterra}. 

The rest of paper is organized as follows.  In Section \ref{sec:trig}, we summarize the relevant results of trigonometric interpolation algorithm developed in \cite{zou_tri}.  Section \ref{sec:tiba_nonlinear} is devoted to develop TIBA, described in Algorithm \ref{alg:nonlinearODE}, for general non-linear system (\ref{eq:nonlinear_ode_order2}-\ref{eq:nonlinear_ode_order2:neum}).  Section \ref{sec:tiba_linear} enhances TIBA with Algorithm \ref{alg:linearODE} for linear second order ODE  (\ref{eq:linear_ode_order2}, \ref{eq:nonlinear_ode_order2:diri}-\ref{eq:nonlinear_ode_order2:neum}),  and some details of derivation are moved to Appendix \ref{appendix:eq:UAV}.  Numerical tests are conducted in Section \ref{sec:performance} to verify established properties of solutions and assess the performance on convergence and accuracy under four sets of boundaries conditions. The summary is made on Section \ref{sec:summary}.

\section{Trigonometric Interpolation on Non-Periodic Functions}\label{sec:trig}
In this section \ref{sec:trig}, we review relevant results of trigonometric interpolation algorithm developed in \cite{zou_tri} starting with following interpolation algorithm on periodic functions. 
\begin{Alg}\label{main_thm} Let $f(x)$ be an odd periodic function \footnote{Similar results for even periodic function is also available in \cite{zou_tri}.} with period $2b$ and $N=2M=2^{q+1}$ for some integer $q\ge 1$ and $x_j, y_j$ are defined by	
	\begin{eqnarray}
		x_j &:=& -b + j\lambda,  \quad \lambda= \frac{2b}N , \quad 0\le j <N, \label{x_grid_N} \\
		y_j &:=& f(x_j),  \label{f_N_interpolation_new} 
	\end{eqnarray}
	then there is a unique $M-1$ degree trigonometric polynomial
	\begin{eqnarray*}
		f_M(x) &=& \sum_{0< j <M}a_j \sin\frac{j\pi x}b, \label{f_M_interpolation_new_odd}\\
		a_j&=&\frac{2}N\sum_{0< k <N} (-1)^j y_k \sin\frac{2\pi j k}{N}, \quad 0< j <M \label{aj_odd}
	\end{eqnarray*}
	such that it fits to all grid points, i.e.
	\[
		f_M(x_{k})=y_{k}, \quad  0\le k <N.
	\]
\end{Alg} 
One can computer coefficients by Inverse Fast Fourier Transform (ifft):
\[	
\{a_j (-1)^j\}_0^{N-1} = 2\times Imag (ifft(\{y_k\}_{k=0}^{N-1})).
\]
Algorithm \ref{main_thm} has been enhanced so it can be applied to a nonperiodic function $f$ whose $K+1$-th derivative $f^{(K+1)}(x)$ exists over a bounded interval $[s,e]$. To seek for a periodic extension with same smoothness, we assume that $f$ can be extended smoothly such that $f^{(K+1)}$ exists and is bounded over $[s-\delta, e+\delta]$ for certain $\delta>0$.  Such a periodic extension of $f$ can be achieved by a cut-off smooth function $h(x)$ with following property: 
\begin{eqnarray*}
	h(x)=\left\{\begin{array}{cc}
		1 & x\in [s,e], \\
		0 & \mbox{$x<s-\delta$ or $x>e+\delta$}. \\
	\end{array}\right.
\end{eqnarray*}
A cut-off function with closed-form analytic expression is proposed  in \cite{zou_tri}.  Let
\begin{equation}\label{ob}
	o=s-\delta, \quad b=e+\delta -o,
\end{equation}
and define $F(x):=h(x+o)f(x+o)$ for $x\in [0,b]$. One can treat $F(x)$ as an odd periodic function with period $2b$. Apply Algorithm \ref{main_thm} to generate the trigonometric interpolation of degree $M-1$ with $N$ evenly-spaced grid points over $[-b,b]$
\[
F_M(x) = \sum_{0 < j < M} a_j \sin\frac{j\pi x}{b},
\]
and let 
\[
	\hat{f}_{M}(x)=F_M(x-o)=\sum_{0< j < M} a_j \sin\frac{j\pi (x-o)}{b}. 
\]
$\hat{f}_{M}(x)|_{[s,e]}$ can be treated as an trigonometric interpolation of $f$ since $\hat{f}_{M}(x_k)=f(x_k)$ for all grid points $x_k\in [s,e]$.  Numerical tests on certain functions demonstrate that $\hat{f}$ approaches to $f$ with decent accuracy \cite{zou_tri}. 
\section{TIBA for Second Order Nonlinear ODE}\label{sec:tiba_nonlinear} To apply trigonometric interpolation for non-periodic functions, we assume that $f$ in Eq (\ref{eq:nonlinear_ode_order2}) is continuous differential on $[s-\delta, e+\delta]\times R^2$ for certain $\delta>0$.  By parallel shifting if needed,  we assume $s=\delta$ without loss of generality. Let $h$  be a cut-off function specified in Section \ref{sec:trig} and 
construct $f_h(x,v,u)$ as follows
\[
f_h(x,v,u) = f(x,v,u)h(x), \qquad (x,v,u) \in [0,b] \times R^2.
\]
Consider a solution $v(x)$ of the following ODE system
\begin{eqnarray}
	v''(x) &=& f_h(x, v,v'),  \quad x\in [0,b] \label{eq:nonlinear_ode_order2_F},\\
	\alpha  &=& d_{11}v(s) + d_{12}v'(s) +d_{13}v(e) + d_{14}v'(e), \label{eq:nonlinear_ode_order2:diri_F}\\ 
	\beta  &=& d_{21}v(s) + d_{22}v'(s) +d_{23}v(e) + d_{24}v'(e). \label{eq:nonlinear_ode_order2:neum_F} 
\end{eqnarray}
It is clear that $v(x)|_{[s,e]}$ solves ODE (\ref{eq:nonlinear_ode_order2}-\ref{eq:nonlinear_ode_order2:neum}). Define $u(x):=v'(x)$ and $z(x):=v''(x)$.  By Eq (\ref{eq:nonlinear_ode_order2_F}),  $z(x)$ and its derivatives $z^{(k)}$ vanish at boundary points $\{0,b \}$, hence it can be smoothly extended as an odd periodic function with period $2b$ and be approximated by trigonometric polynomial.  Assume 
\begin{equation}\label{eq:z_M_nonlinear}
	\tilde{z}_M(x) = \sum_{0\le j<M}b_j \sin \frac{j\pi x}{b} 
\end{equation}
is an interpolant of $z(x)$ with $N$ equispaced grid points over $[-b,b]$ by Algorithm \ref{main_thm}. $u$ and $v$ can be derived accordingly 
\begin{eqnarray}
	\tilde{u}_M(x ) &=& a_0 -   \frac{b}{\pi }\sum_{1\le j <M} \frac{b_j}{j} \cos  \frac{j\pi x}{b},  \label{eq:tildeu}\\
	\tilde{v}_M(x ) &=& a_1 +  a_0 x -  (\frac{b}{\pi })^2 \sum_{1\le j <M} \frac{b_j}{j^2}  \sin  \frac{j\pi x}{b},  \label{eq:tildev_nonlinear}
\end{eqnarray}
where $a_0,a_1$ are two constant and can be determined by boundary conditions as shown in Eq (\ref{eq:a01}) below. 

The following notations and conventions will be adopted in the rest of this paper. A $k$-dim vector is considered as $(k,1)$ dimensional matrix unless specified otherwise. Define
\begin{eqnarray*} 
	x_k &=& k \lambda, \quad \lambda = \frac{b}M ,\quad 0\le k \le M,  \quad X = (x_k)_{0\le k\le M},\\
	u_k &=& \tilde{u}_M(x_k), \quad v_k = \tilde{v}_M(x_k),  \quad z_k = \tilde{z}_M(x_k),  \quad f_k = f_h(x_k, v_k, u_k),\\
	U &=& (u_k)_{0\le k\le M},\quad V = (v_k)_{0\le k\le M}, \quad Z = (z_k)_{0\le k\le M}, \\
	F &=& (f_k)_{0\le k\le M} , \quad K = (1,2,\cdots, M-1)^T, \quad B=(b_i)_{1\le i<M},\\
	I &=& (1,1,\cdots, 1)^T_{M-1}, \quad I_a = (-1,1, -1,\cdots, -1)^T_{M-1}.
\end{eqnarray*}
For any two matrices $A,B$ with same shape,  $A\circ B$ denotes the Hadamard product, which applies the element-wise product to two matrices. $AB$ denote the standard matrix multiplication when applicable. $\sum(W)$ denotes the sum of all elements in a vector $W$. $A(i,:)$ and $A(:,j)$ is used to denote the $i$-th row and $j$-th column of $A$ respectively.  $W(k:l)$ denote the $l-k+1$-th vector $(w_k, \cdots, w_l)^T$. In addition, $diag (W)$ is the diagonal matrix constructed by $W$.  Note we have
\[
s = x_{m}, \qquad e=x_{m+n}.
\]
At grid points of interpolation, ODE dynamic (\ref{eq:nonlinear_ode_order2_F}) is characterized by
\begin{equation} \label{eq:bdp_orde2_d_nonlinear} 
	Z = F.  
\end{equation}
$Z, U, V$ can be calculated based on Eq. (\ref{eq:z_M_nonlinear}-\ref{eq:tildev_nonlinear}):
\begin{eqnarray}
	z_k &=& \sum_{0\le j<M}b_j \sin \frac{2\pi jk}{N},   \label{eq:z_M_d}\\
	u_k &=& a_0 - \frac{b}{\pi }\sum_{1\le j <M} \frac{b_j}{j} \cos  \frac{2\pi jk }{N},  \label{eq:tildeu_d}\\
	v_k &=& a_1 + a_0 x_k-(\frac{b}{\pi })^2 \sum_{1\le j <M}\frac{b_j}{j^2}\sin \frac{2\pi jk }{N}. \label{eq:tildev_d}
\end{eqnarray}
Eq (\ref{eq:tildev_d}) can be used to solve $a_0$ and $a_1$:
\begin{equation}\label{eq:a01}
	a_0 = \frac{v_M - v_0}b, \quad a_1 = v_0.
\end{equation}
Define 
\begin{equation}\label{eq:SC}
	S =(\sin\frac{2\pi jk}{N})_{1\le j,k <M}, \qquad C =  (\cos\frac{2\pi jk}{N})_{1\le j,k < M} .
\end{equation}
It is easy to check $SS = \frac{M}2 E$, where $E$ is the $M-1$ identity matrix, and therefore $O:=\sqrt{\frac M2}S$ is a symmetric orthogonal matrix. Define
\[
\Theta = (\theta_{ij})_{1\le i, j < M} = O \cdot diag(1/K^2) \cdot O.
\]
We need represent $B, U$ in term of $V$. First, rewrite (\ref{eq:tildev_d}) in vector format:
\begin{equation}\label{eq:v_vector}
	V(1:M-1) = a_1I +a_0 \frac bM K - (\frac{b}{\pi })^2 S \cdot diag(1/K^2)\cdot B,
\end{equation}
which implies
\begin{equation}\label{eq:B}
	B = diag(K^2) S (\frac{2a_1 \pi^2}{Mb^2} I + \frac{2a_0 \pi^2}{bM^2}K - \frac{2\pi^2}{Mb^2} V(1:M-1)).
\end{equation}
Applying Eq (\ref{eq:B}) and Eq (\ref{eq:a01}) to (\ref{eq:z_M_d}), we obtain
	\[
	Z(1:M-1) = \frac{v_0\pi^2}{Mb^2}\Theta^{-1} (MI-K) + \frac{v_M\pi^2}{M b^2}\Theta^{-1} K -\frac{\pi^2}{b^2}\Theta^{-1} V(1:M-1),
	\]
and a discretization of Eq (\ref{eq:bdp_orde2_d_nonlinear})
\begin{equation}\label{eq:ode_discrete_nonlinear}
	\frac{v_0\pi^2}{Mb^2} (MI-K) + \frac{v_M\pi^2}{M b^2} K -\frac{\pi^2}{b^2} V(1:M-1) - \Theta \cdot F(1:M-1) = 0.
\end{equation}
By Eq (\ref{eq:tildeu_d}), we can represent $U$ by linear combination of $V$ 
\begin{equation}\label{eq:UAV}
	U = AV.
\end{equation}
We leave details of derivation of $A=(a_{ij})_{0\le i, j\le M}$ in Appendix \ref{appendix:eq:UAV} and show the results as follows.
\begin{eqnarray}
	a_{0,0} &=& \frac{\pi}{b} sum(I_a \circ \cot(\pi K/N))\nonumber\\
	&-&\frac{\pi}{bM} sum(I_a\circ K\circ \cot(\pi K/N))-\frac1b \label{eq:top}\\
	a_{0,1:M-1} &=& -\frac{\pi}{b} I'_a \circ \cot(\pi K'/N) \nonumber \\
	a_{0,M} &=& \frac{\pi}{bM} sum(I_a\circ K\circ \cot(\pi K/N))+\frac1b \nonumber\\
	a_{i,0} &=& \frac{\pi}{2b} sum((-1)^i \cot(i,:)I_a) \nonumber\\
	&-& \frac{\pi}{2bM} sum( (-1)^i I_a \circ \cot(i,:) \circ K) -1/b, \label{eq:middle}\\
	a_{i,1:M-1} &=& \frac{\pi}{2b}(-1)^{i+1} I_a' \circ \cot(i,:), \nonumber \\
	a_{i,M} &=& \frac{\pi}{2bM} sum((-1)^i I_a \circ \cot(i,:) \circ K)+1/b. \nonumber \\
	a_{M,0} &=& -\frac{\pi}{b} sum(I_a \circ \tan(\pi K/N)) \nonumber\\
	&-&\frac{\pi}{bM} sum((K\circ \cot(\pi K/N))-\frac1b, \label{eq:bottom}\\
	a_{M,1:M-1} &=& \frac{\pi}{b} I'_a \circ \tan(\pi K'/N), \nonumber\\
	a_{M,M} &=& \frac{\pi}{bM} sum((K\circ \cot(\pi K/N))+\frac1b,  \nonumber
\end{eqnarray}
where $\quad 0<i<M$ and
\[
\cot(k, i) := Cot\frac{k+i}{N}\pi + Cot\frac{k-i}{N}\pi,
\]
and $Cot(x)=\cot(x)$ if $x/\pi$ is not integer and $Cot(x)=0 $ otherwise. 

Note that $A$ is singular since $\sum_{j}a_{i,j}=0$ for all $0\le i\le M$, and $V$ can not be recovered by $U$, which is expected. 
By Eq (\ref{eq:UAV}),  we obtain a non-linear system on $V$ by combining boundary conditions (\ref{eq:nonlinear_ode_order2:diri_F}-\ref{eq:nonlinear_ode_order2:neum_F}) and the discretization  (\ref{eq:ode_discrete_nonlinear}) of Eq  (\ref{eq:nonlinear_ode_order2_F}):
\begin{eqnarray}
\quad h_0(V) &:=&d_{11}v_m + d_{13}v_{n+m} \nonumber\\
&+& \sum_{0\le i\le M} (d_{12}a_{m,k} +d_{14}a_{m+n,k})v_k -\alpha =0  \label{eq:alpha}\\
\quad h_i(V)&:=& \frac{v_0\pi^2}{Mb^2}(M-i) + \frac{v_M\pi^2}{M b^2}i -\frac{\pi^2}{b^2} v_i \nonumber\\
&-& \sum_{1\le j <M} \theta_{ij} f_j = 0, \quad 0<i<M,  \label{eq:nonlinear_middle}\\
\quad h_M(W)&:=&d_{21}v_m + d_{23}v_{n+m} \nonumber\\
&+& \sum_{0\le i\le M} (d_{22}a_{m,k} +d_{24}a_{m+n,k})v_k - \beta = 0,  \label{eq:beta}
\end{eqnarray}
with
\[
f_j =h(x_j)f(x_j, v_j, u_j),\qquad u_j =\sum_{0\le k\le M} a_{jk} v_k.
\]
TIBA can be summarized as follows
\begin{Alg}\label{alg:nonlinearODE}
\begin{enumerate}
	\item Select proper $N, \delta$ and follow Section \ref{sec:trig} to generate cut-off function $h$; 
	\item Apply Eq (\ref{eq:top}-\ref{eq:bottom}) to generate matrix $A$;
	\item Solve $M+1$-dim system that consists of $M+1$ equations (\ref{eq:alpha}-\ref{eq:beta}).
	\item Apply Eq (\ref{eq:a01}) and (\ref{eq:B}) to calculate parameters $a_0, a_1, b_1, \dots, b_{M-1}$ of $\tilde{v}_M$ defined by Eq. (\ref{eq:tildev_nonlinear}); 
	\item Restrict $\tilde{v}_M$ to $[s,e]$ as approximation of target solution. 
\end{enumerate}
\end{Alg}
In general, the performance of Algorithm \ref{alg:nonlinearODE} depends on how effective 
the nonlinear system (\ref{eq:alpha}-\ref{eq:beta}) can be solved.  Note that close form of Jacobian $\frac{\partial h_i}{\partial v_j}$ is available and classic Newton method can be applied. TIBA becomes particular attractive when ODE (\ref{eq:nonlinear_ode_order2}) is reduced to a linear system as discussed in Section \ref{sec:tiba_linear}.
\section{TIBA for Second Order Linear ODE}\label{sec:tiba_linear}
In this section,  we focus on second order linear ODE and assume $f(x,v,u)$ in Eq. (\ref{eq:nonlinear_ode_order2}) is a linear function in $u,v$ as follows
\[
f(x,v,u) = p(x)u+q(x)v(x) + r(x),
\]
and $p,q,r$ is continuous differential on $[s-\delta, e+\delta]$.  Eq. (\ref{eq:nonlinear_ode_order2}) is reduced to
\begin{equation}\label{eq:linear_ode_order2}  
	y''(x) = p(x)y'+q(x)y(x) +r(x), \quad x\in [s,e].
\end{equation}
We follow same conventions and notations as in Section \ref{sec:tiba_nonlinear}. In addition, define 
\[
p_h(x)=p(x)h(x), \quad q_h(x)=q(x)h(x), \quad  r_h(x) = r(x)h(x),  \quad x \in [0,b]
\]
and
\begin{eqnarray*} 
	p_{k} &=& p_h(x_k), \quad q_k = q_h(x_k), \quad r_k= r_h(x_k),\\
	P&=&(p_k)_{0< k < M}, \quad  Q= (q_k)_{0< k< M},\quad R = (r_k)_{0<k< M},\\
	\hat{U} &=& U(1:M-1), \quad \hat{V} = V(1:M-1).
\end{eqnarray*}
Replacing $F(1:M-1)$ by $P\circ U + Q\circ V + R$ in Eq (\ref{eq:ode_discrete_nonlinear}) and applying (\ref{eq:UAV}), we obtain 
\begin{equation}\label{eq:ode_discrete}
	\frac{v_0\pi^2}{Mb^2} (MI-K) + \frac{v_M\pi^2}{M b^2} K -\frac{\pi^2}{b^2} \hat{V} = \Theta \cdot (R + Q\circ \hat{V}  + P\circ \hat{U})
\end{equation}
Eq (\ref{eq:nonlinear_ode_order2:diri_F}-\ref{eq:nonlinear_ode_order2:neum_F}) can be written as two linear equations in $V$. Insert Eq (\ref{eq:nonlinear_ode_order2:diri_F}) in front and attach \ref{eq:nonlinear_ode_order2:neum_F} to end of linear system (\ref{eq:ode_discrete}), we obtain a $M+1$-dim linear system 
\begin{equation}\label{eq:PhiVG}
	\Phi V = \Psi,
\end{equation}
where $\Psi=(\psi_i)_{0\le i\le M}$ is determined by
\begin{equation} \label{eq:psi}
	\psi_0=\alpha, \quad \psi_M= \beta, \quad  \Psi(1:M-1)= -\Theta R;
\end{equation}
$\Phi=(\phi_{ij})_{0\le i,k \le M}$ is determined by ($ 0\le k \le M $),
\begin{eqnarray}
	\phi_{0,k} &=& d_{12}\cdot A(m,k) + d_{14}\cdot A(m+n,k) \nonumber\\
	&+& \delta_{m,k}d_{11} + \delta_{m+n,k}d_{13}, \label{eq:phi_0} \\
	\phi_{M,k} &=& d_{22}\cdot A(m,k) + d_{24}\cdot A(m+n,k) \nonumber\\
	&+& \delta_{m,k}d_{21} + \delta_{m+n,k}d_{23}, \label{eq:phi_M} 
\end{eqnarray}
and for $0< i <M$, 
\begin{eqnarray}
	\phi_{i,0} &=& -\frac{(M-i)\pi^2}{Mb^2} + (\theta(i,:)\circ P) \cdot A(1:M-1,0),\label{eq:phi_i_0} \\
	\phi_{i,M} &=& -\frac{i \pi^2}{b^2 M} + (\theta(i,:) \circ P) \cdot A(1:M-1, M),\label{eq:phi_i_M} \\
	\phi(i,1:M-1) &=& \frac{\pi^2}{b^2}(\delta_{ij})_{1\le j<M} +\theta(i,:)\circ Q^T  \nonumber\\
	&+&(\theta(i,:) \circ P^T) \cdot A(1:M-1,1:M-1) \label{Phi}  \label{eq:phi_i_middle}
\end{eqnarray}
For any solution $y$ of ODE system (\ref{eq:linear_ode_order2},\ref{eq:nonlinear_ode_order2:diri}-\ref{eq:nonlinear_ode_order2:neum}), 
Eq (\ref{eq:PhiVG}) always holds .  As such, we have
\begin{Thm} \label{theorem_uniquness}
	Let $\Phi$ be defined by Eq (\ref{eq:phi_0}-\ref{eq:phi_i_middle}) and assume that there are solutions for ODE  (\ref{eq:linear_ode_order2}) with boundary conditions (\ref{eq:nonlinear_ode_order2:diri}-\ref{eq:nonlinear_ode_order2:neum}), then solution is unique if $rank(\Phi)=M+1$. Furthermore, there is no solution if $M+1>rank(\Phi, g)>rank(\Phi)$.
\end{Thm}
\begin{Rem} 
If $rank(\Phi)=M+1$, the identified $\tilde{v}_M$ by Eq. (\ref{eq:tildev_nonlinear}) is expected to converge a solution based on Algorithm \ref{main_thm}. In addition, it is not hard to verify numerically whether $\tilde{v}_M$ meets Eq. (\ref{eq:linear_ode_order2}) as shown in Section \ref{sec:performance}. Note that boundary conditions (\ref{eq:nonlinear_ode_order2:diri},\ref{eq:nonlinear_ode_order2:neum}) are always satisfied by $\tilde{v}_M(x)$. Similarly, it is expected that there are infinite many solutions if $rank(\Phi, g)=rank(\Phi)<M+1$ since there are infinitely many $\tilde{v}_M(x)$. We provides numerical tests to demonstrate possible scenarios of solution structure in Section \ref{sec:performance}.
\end{Rem}
 Algorithm \ref{alg:nonlinearODE} can be updated to solve ODE  (\ref{eq:linear_ode_order2},\ref{eq:nonlinear_ode_order2:diri}-\ref{eq:nonlinear_ode_order2:neum}) as follows.
\begin{Alg}\label{alg:linearODE}
	\begin{enumerate}
		\item Follow step 1 and 2 in Algorithm \ref{alg:nonlinearODE};
		\item Calculate $\Phi, \Psi$ by Eq (\ref{eq:psi}) and (\ref{eq:phi_0}-\ref{eq:phi_i_middle}) and further obtain $V$ by Eq. (\ref{eq:PhiVG});
		\item Follow step 4, 5 in Algorithm \ref{alg:nonlinearODE};
	\end{enumerate}
\end{Alg}

\section{Performance}\label{sec:performance}
The tests in this section include four sets of boundary conditions shown in Table \ref{tab:test}. 
\begin{table}[htbp]
	\tiny
	\caption{The types of boundary conditions. $\{d_{ij}\}_{1\le i,j\le 4}$ are parameters in Eq (\ref{eq:nonlinear_ode_order2:diri}-\ref{eq:nonlinear_ode_order2:neum}).}
	\begin{tabular}{lrrrrrrrrl}
		type	& \multicolumn{1}{l}{$d_{11}$} & \multicolumn{1}{l}{$d_{12}$} & \multicolumn{1}{l}{$d_{13}$} & \multicolumn{1}{l}{$d_{14}$} & \multicolumn{1}{l}{$d_{21}$} & \multicolumn{1}{l}{$d_{22}$} & \multicolumn{1}{l}{$d_{23}$} & \multicolumn{1}{l}{$d_{24}$} & condition on \\ \hline\hline
		$Neumann$ & 1     & 0     & 0     & 0     & 0     & 1     & 0     & 0     & $v_s,u_s$ \\
		$Dirichlet$ & 1     & 0     & 0     & 0     & 0     & 0     & 1     & 0     & $v_s, v_e$ \\
		$Mix_1$   & 1     & 0     & 0     & 0     & 0     & 0     & 0     & 1     & $v_s, u_e$ \\ 
		$Mix_2$ & 1     & 1     & 0     & 0     & 0     & 0     & 1     & 1     & $v_s+u_s, v_e+u_e$ \\ \hline
	\end{tabular}%
	\label{tab:test}%
\end{table}%
We use both TIBO and classic Runge-Kutta scheme (Section 9.4.1, page 284 \cite{ryts}), labeled $rk4$, as benchmarks for performance assessment on $Neumann$ type.  For $Dirichlet$ and $Mix_1$, we follow the shooting method, i.e. searching for $y'(s)$ to meet required boundary conditions; then apply Runge-Kutta scheme. rk4's result is further used as the initial values for the optimization process used in TIBO.  Same benchmarks are applied for $mix_2$ except we need search both $y(s)$ and $y'(s)$ to meet Equations  (\ref{eq:nonlinear_ode_order2:diri}-\ref{eq:nonlinear_ode_order2:neum}).
The performance of shooting method is sensitive to the initial guess used in the shooting method.  For each of four test types, we conduct $5$ tests with randomly selected initial values $y(s),y'(s)$ 
and adjust $init_y$ to the given $y(s)$ for type $Dirichlet$ and $Mix_1$.
\subsection{Test on existence and uniqueness}\label{subsec:uniquness}
In this section,  we conduct numerical tests to verify the property of ODE system (\ref{eq:linear_ode_order2},\ref{eq:nonlinear_ode_order2:diri}-\ref{eq:nonlinear_ode_order2:neum}) by Theorem \ref{theorem_uniquness}.  It is shown in \cite{ode_naga} that the existence and uniqueness issue of BVP for a non-homogeneous linear ODE is the same as that for the associated homogeneous ODE. The tests in this subsection will be based on following homogeneous linear ODE system:
\begin{eqnarray}
	y''(x) +2\pi y' + \frac54 \pi^2 y &=& 0 \label{eq:example_uniquness} \\
 d_{11}y(s) + d_{12}y'(s) +d_{13}y(e) + d_{14}y'(e) &=&\alpha, \label{eq:uniquness_cond1}\\ 
d_{21}y(s) + d_{22}y'(s) +d_{23}y(e) + d_{24}y'(e) &=& \beta. \label{eq:uniquness_cond2} 
\end{eqnarray}
One can solve Eq (\ref{eq:example_uniquness}) and obtain following two special solutions.
\begin{eqnarray*}
	y_1(x) &=& e^{-(x-1)\pi}(\cos\frac{\pi(x-1)}{2} + 2\sin\frac{\pi(x-1)}{2})\\
	y_2(x) &=& e^{-(x-1)\pi} \sin\frac{\pi(x-1)}{2},
\end{eqnarray*}
where $s=1, e=3$. 
Other solution of Eq (\ref{eq:example_uniquness}) can be expressed by
\[
y(x) = c_1y_1(x)+c_2y_2(x), 
\]
and relevant boundary values are
\begin{eqnarray*}
y(s)&=&c_1,\quad y'(s) = \frac{\pi}2c_2, \quad y(e) = -e^{-2\pi}c_1, \quad y_1'(e)=-\frac{\pi}2 e^{-2\pi}c_2.
\end{eqnarray*}
Therefore it is expected that
\begin{enumerate}
	\item on type $Neumann$, there is a unique solution with $c_1 = \alpha$ and $c_2=\frac{\pi}2\beta$; 
	\item on type $Dirichlet$,  there are infinitely many solution if $\frac{\beta}{\alpha}=-e^{2\pi}$ with $c_1=\alpha$ and any $c_2$; there is no solution if $\frac{\beta}{\alpha}\neq -e^{2\pi}$;
	\item on type $Mix_1$,  there is a unique solution with $c_1 = \alpha$ and $c_2=-\frac2{\pi} e^{2\pi}\beta$;
	\item on type $Mix_2$,  there are infinitely many solution if $\frac{\beta}{\alpha}=-e^{2\pi}$ with any $c_1, c_2$ such that $c_1+c_2\pi/2=\alpha$; there is no solution if $\frac{\beta}{\alpha}\neq -e^{2\pi}$. 
\end{enumerate}
Define $y_b(x)=y_1(x)+y_2(x)$, and let $\alpha_b,\beta_b$ be derived by boundary conditions (\ref{eq:uniquness_cond1}-\ref{eq:uniquness_cond2}), i.e.
\begin{eqnarray*}
\alpha_b  &=& d_{11}y_b(s) + d_{12}y_b'(s) +d_{13}y_b(e) + d_{14}y_b'(e), \label{eq:uniquness_cond1_b}\\
\beta_b  &=& d_{21}y_b(s) + d_{22}y_b'(s) +d_{23}y_b(e) + d_{24}y_b'(e). \label{eq:uniquness_cond2_b}
\end{eqnarray*}
We call $y_b$ as the base solution of ODE (\ref{eq:example_uniquness}) with above boundary conditions.

We report some of following measures for each of four test types in Table \ref{tab:test}. 
\begin{enumerate}
	\item $max|y_o''-f|$:  The max difference is calculated by TIBO-based solution $y_{o}$, i.e. $max|y_{o}''-f(x,y_{o}, y'_{o})|$.  The max is taken over the set of $M=2^{10}$ equally-spaced $x$ values over $[0,b]$.  Note that grid point set used in the optimization algorithm is determined by $M=2^{7}$. As such, negligible $max|y_o''-f|$ indicates that $y_o$ converges to the solution over $[0,b]$.
	\item  $max|y_{l}-y_{b}|$,  $max|y_{o}-y_{b}|$ and  $max|y_{rk4}-y_{b}|$:  the max differences between base $y_b$ and  three approximations by  TIBA, TIBO, and rk4 respectively.  The max is taken over all grid points over $[s,e]$ determined by $M=2^7$.
\end{enumerate}
\subsubsection{Type $Neumann$}\label{subsec:Neumann}
Table \ref{tab:linear_test_neum} summarizes the test results with $\alpha=\alpha_b, \beta=\beta_b$.  Both $y_{l}$ and $y_{o}$ converges to the unique base solution $y_b$ with almost same accuracy and they outperform $y_{rk4}$ significantly.  The negligible  difference between $y_{l}$ and $y_{o}$ suggests that TIBO converges effectively, which is also confirmed by negligible value of $max|y_o''-f|$.
\begin{table}[htbp]
	\small
	\caption{The performance on $Neumann$ condition. }
	\begin{tabular}{lccc}
		\multicolumn{1}{l}{$max|y_{l}-y_{b}|$} & \multicolumn{1}{l}{$max|y_{o}-y_{b}|$} & \multicolumn{1}{l}{$max|y_{rk4}-y_{b}|$} & \multicolumn{1}{l}{$max|y_o''-f|$} \\ \hline \hline
		1.5E-09 & 1.7E-09 & 1.7E-06 & 1.2E-06 \\ \hline
	\end{tabular}%
	\label{tab:linear_test_neum}%
\end{table}%
\subsubsection{Type $Dirichlet$}\label{subsec:Dirichlet}
Table \ref{tab:linear_test_diri} summarizes the test results with $\alpha=\alpha_b, \beta=\beta_b$ and therefore condition $\frac{\beta}{\alpha}=-e^{2\pi}$ is met and infinite many solutions are expected.  Negligible $max|y_o''-f|$ and different $max|y_{o}-y_{b}|$ show that $y_o$ converges to different values on each of $5$ tests, consistent to what is expected. $y_l$ is not dependable as shown by significant $max|y_{l}-y_{b}|$.  $\Phi$ in Eq (\ref{eq:PhiVG}) turns out to be very close to a singular matrix by a {\em matlab} warning message ``Matrix is close to singular or badly scaled". 

To test other situation where $\frac{\beta}{\alpha}\neq -e^{2\pi}$ is not satisfied,  we keep $\alpha=\alpha_b$ and replace $\beta$ by $1.1\beta_b$. Table \ref{tab:linear_test_diri} summarizes the results.  The value of $y_l$ is clearly out of range and $max|y_o''-f|$ become significantly bigger than before, suggesting that there is no solution. 
\begin{table}[htbp]
	\small
	\caption{The performance on $Dirichlet$ condition. Case 1: Using conditions (\ref{eq:uniquness_cond1_b}-\ref{eq:uniquness_cond2_b}). Infinite many solutions due to $\frac{\beta_b}{\alpha_a}=-e^{2\pi}$.  Five tests are conducted based on the combination of shooting and RK5 with randomly selected  initial value of $y'(s)$.}
	\begin{tabular}{lccc}
		\multicolumn{1}{l}{Test Id} & \multicolumn{1}{l}{$max|y_{l}-y_{b}|$} & \multicolumn{1}{l}{$max|y_{o}-y_{b}|$}  & \multicolumn{1}{l}{$max|y_o''-f|$} \\ \hline \hline
		1     & 4.0E-01 & 3.5E-02 &  1.3E-06 \\
		2     & 4.0E-01 & 4.2E-02 &  1.1E-06 \\
		3     & 4.0E-01 & 1.5E-02 &  1.3E-06 \\
		4     & 4.0E-01 & 2.5E-02 &  1.2E-06 \\
		5     & 4.0E-01 & 5.1E-02 &  1.4E-06 \\ \hline 
	\end{tabular}%
	\label{tab:linear_test_diri}%
\end{table}%
\begin{table}[htbp]
	\small
	\caption{The performance on $Dirichlet$ condition. Case 2: keep $\alpha_a$ and replace $\beta_b$ by $1.1\beta_b$ in conditions (\ref{eq:uniquness_cond1_b}-\ref{eq:uniquness_cond2_b}). Five tests are conducted based on same method as in Case 1}
	\begin{tabular}{lccc}
		\multicolumn{1}{l}{Test Id} & \multicolumn{1}{l}{$max|y_{l}-y_{b}|$} & \multicolumn{1}{l}{$max|y_{o}-y_{b}|$}  & \multicolumn{1}{l}{$max|y_o''-f|$} \\ \hline \hline
		1     & 2.3E+08 & 4.0E-02 & 3.3E-03 \\
		2     & 2.3E+08 & 3.8E-02 & 3.3E-03 \\
		3     & 2.3E+08 & 2.0E-02 & 3.3E-03 \\
		4     & 2.3E+08 & 2.1E-02 & 3.3E-03 \\
		5     & 2.3E+08 & 5.7E-02 & 3.3E-03 \\ \hline 
	\end{tabular}%
	\label{tab:linear_test_diri_no_solution}%
\end{table}%
\subsubsection{Type $Mix_1$}\label{subsec:mix_1}
Table \ref{tab:linear_test_mix_1_wo_neum_init} summarizes the test results with $\alpha=\alpha_b, \beta=\beta_b$. 
We add $|y'_o(s)-y'_b(s)|$ and $|y_o(e)-y_b(e)|$ in the table to show that shooting method doesn't generate an accurate estimation on $y'(s)$, which deteriorates TIBO's performance with visible error $max|y_{o}-y_{b}|$.

TIBA provides a decent approximation of the unique solution $y_b$ as shown by negligible $max|y_{l}-y_{b}|$ and outperforms TIBO quite dramatically.  

To enhance the performance of TIBO,  we improve the quality of initial values consumed by the optimization process following three steps: first, use shooting method to estimate $y'(s)$; secondly, apply TIBO with $Neumann$ condition with given $y(s)$ and estimated $y'(s)$, finally take the output of second step as initial values of TIBO  with $Mix_1$ condition.  The result is shown in Table \ref{tab:linear_test_mix_2}. The performance of TIBO has been significantly improved as shown by $max|y_{o}-y_{b}|$,  $|y_o(e)-y_b(e)|$ and $|y'_o(s)-y'_b(s)|$.  Nevertheless, TIBA still outperforms TIBO significantly. 
\begin{table}[htbp]
	\small
	\caption{The performance on $Mix_1$ condition: TIBA vs TIBO with shooting-rk4 based initial values}
	\begin{tabular}{lcccc}
		\multicolumn{1}{l}{ID} & \multicolumn{1}{l}{$max|y_{l}-y_{b}|$} & \multicolumn{1}{l}{$max|y_{o}-y_{b}|$} &  \multicolumn{1}{l}{$max|y_o''-f|$} & \multicolumn{1}{l}{$|y'_o(s)-y'_b(s)|$}  \\ \hline \hline
		1     & 4.1E-08 & 4.6E-03 &  9.2E-04 & 4.1E-02  \\ \hline
	\end{tabular}%
	\label{tab:linear_test_mix_1_wo_neum_init}%
\end{table}%
\begin{table}[htbp]
	\small
	\caption{The performance on $Mix_1$ condition: TIBA vs TIBO with enhanced initial values (Subsection \ref{subsec:mix_1})}
	\begin{tabular}{lcccc}
		\multicolumn{1}{l}{ID} & \multicolumn{1}{l}{$max(|y_{l}-y_{b}|)$} & \multicolumn{1}{l}{$max(|y_{o}-y_{b}|)$}  & \multicolumn{1}{l}{$max|y_o''-f|$} & \multicolumn{1}{l}{$|y'_o(s)-y'_b(s)|$} \\ \hline \hline
		1     & 4.1E-08 & 6.9E-06 & 1.4E-06 & 6.1E-05   \\ \hline
	\end{tabular}%
	\label{tab:linear_test_mix_1_w_neum_init}%
\end{table}%
\subsubsection{Type $Mix_2$}
The distribution of solutions conditioning on $Mix_2$ is similar on Dirichlet.  Tests are conducted to test two cases: infinite many solutions vs no solution.  The results are shown in Table \ref{tab:linear_test_mix_2} and \ref{tab:linear_test_mix_2_no_solution} respectively.  In both cases, TIBA returns no solution as confirmed by $max|y_l-y_b|$. TIBO returns five different solutions with $\alpha=\alpha_b, \beta=\beta_b$ as indicated by negligible $max|y_o''-f|$ in Table \ref{tab:linear_test_mix_2}  and no solution with $\alpha=\alpha_b, \beta=1.1\beta_b$ as implied by noticeable $max|y_o''-f|$ in Table \ref{tab:linear_test_mix_2_no_solution}.
\begin{table}[htbp]
	\small
	\caption{The performance on $Mix_2$ condition. TIBA vs TIBO with  $\alpha=\alpha_b, \beta=\beta_b$.}
	\begin{tabular}{rrrr}
		\multicolumn{1}{l}{Test ID} & \multicolumn{1}{l}{$max(|y_{l}-y_{b}|)$} & \multicolumn{1}{l}{$max(|y_{o}-y_{b}|)$} & \multicolumn{1}{l}{$max|y_o''-f|$} \\ \hline \hline
		1     & 3.5E+00 & 3.2E-01 & 1.1E-06 \\
		2     & 3.5E+00 & 3.7E-01 & 1.4E-06 \\
		3     & 3.5E+00 & 1.3E-01 & 1.2E-06 \\
		4     & 3.5E+00 & 2.2E-01 & 1.3E-06 \\
		5     & 3.5E+00 & 4.6E-01 & 1.0E-06 \\ \hline
	\end{tabular}%
	\label{tab:linear_test_mix_2}%
\end{table}%

\begin{table}[htbp]
	\small
	\caption{The performance on $Mix_2$ condition. TIBA vs TIBO with  $\alpha=\alpha_b, \beta=1.1\beta_b$.}
	\begin{tabular}{rrrr}
		\multicolumn{1}{l}{Test ID} & \multicolumn{1}{l}{$max(|y_{l}-y_{b}|)$} & \multicolumn{1}{l}{$max(|y_{o}-y_{b}|)$} & \multicolumn{1}{l}{$max|y_o''-f|$} \\ \hline \hline
		1     & 2.4E+06 & 3.0E-01 & 5.5E-03 \\
		2     & 2.4E+06 & 4.0E-01 & 5.5E-03 \\
		3     & 2.4E+06 & 1.1E-01 & 5.5E-03 \\
		4     & 2.4E+06 & 2.5E-01 & 5.5E-03 \\
		5     & 2.4E+06 & 4.4E-01 & 5.5E-03 \\ \hline
	\end{tabular}%
	\label{tab:linear_test_mix_2_no_solution}%
\end{table}%
\newpage
\subsection{Test on convergence}
We conduct convergence tests for $Nenumann$ and $Mix_1$ whose base $y_b$ is the unique solution and show the results in Table \ref{tab:linear_test_neum_convergence} and \ref{tab:linear_test_mix_1_convergence} respectively.  For $Nenumann$, we present errors of TIBA, TIBO and rk4. 
For $Mix_1$,  we apply the enhanced initial values in TIBO as described in Subsection \ref{subsec:mix_1} and report the results with Test ID $1$.  One can observe
\begin{enumerate}
	\item TIBA and TIBO converge in a similar pattern and decent accuracy can be achieved starting with $q=7$.  
	\item TIBA outperforms significantly TIBO  for $Mix_1$ type for all cases.  For Nenumann,  TIBA starts to outperform TIBO when $q$ reaches to $8$, suggesting that optimization noise can have non-negligible impact on TIBO's performance when  number of parameters reach to certain level.  
	On the other hand,  TIBA performance is only limited to machine's limitation to handle round-off error.   
	\item TIBA converges much faster than rk4 for $Nenumann$ type,  TIBO also outperforms rk4 although it losses momentum when $q$ reaches to $8$.  
\end{enumerate}
\begin{table}[htbp]
	\small
	\caption{Convergence test on $Nenumann$}
	\begin{tabular}{rrrr}
		\multicolumn{1}{l}{q} & \multicolumn{1}{l}{$max|y_{l}-y_{b}|$} & \multicolumn{1}{l}{$max|y_{o}-y_{b}|$} & \multicolumn{1}{l}{$max|y_{rk4}-y_{b}|$}  \\ \hline\hline
		6     & 1.8E-06 & 1.8E-06 & 3.0E-05 \\
		7     & 1.5E-09 & 1.7E-09 & 1.7E-06  \\
		8     & 1.6E-12 & 3.9E-10 & 1.0E-07  \\
		9     & 1.3E-12 & 9.8E-10 & 6.3E-09  \\ \hline
	\end{tabular}%
	\label{tab:linear_test_neum_convergence}%
\end{table}%
\begin{table}[htbp]
	\small
	\caption{Convergence test on $Mix_1$}
	\begin{tabular}{rrr}
		\multicolumn{1}{l}{q} & \multicolumn{1}{l}{$max|y_{l}-y_{b}|$} & \multicolumn{1}{l}{$max|y_{o}-y_{b}|$} \\ \hline\hline
		6     & 7.7E-05 & 1.1E-04  \\
		7     & 4.1E-08 & 6.9E-06  \\
		8     & 1.2E-10 & 4.3E-07 \\
		9     & 8.0E-11 & 2.7E-08  \\ \hline
	\end{tabular}%
	\label{tab:linear_test_mix_1_convergence}%
\end{table}%
\newpage
\subsection{Test on non-homogeneous functions}
For any given function $f(x)$ over $[s,e]$, one can construct a ODE 
\begin{eqnarray*}
	y'' &=& f''(x)  - f'(x) - f(x) +  p(x) y' +  q(x)y \\
	\alpha  &=& d_{11}f(s) + d_{12}f'(s) + d_{13}f(e) + d_{14}f'(e) \\
\beta  &=& d_{21}f(s) + d_{22}f'(s) + d_{23}f(e) + d_{24}f'(e) 
\end{eqnarray*}
such that $f$ is a solution.  In this section, we test non-homogeneous ODEs with the following base solution used in \cite{zou_tri_III} 
\[
f(x; \theta) =x\cos\theta x, \qquad x\in [1, 3].
\]
with $\theta\in (\pi/2, 3\pi/2)$.  As in \cite{zou_tri_III},  constant $p=0.1$ and $q=1$ are applied.  
TIBO generates close results for five scenarios under each pair $(type, \theta)$ of boundary condition types and function parameters.  Table \ref{tab:linear_non_homogeneous} shows the results with test ID $1$.  Note that $y_{rk4}$ refers to the solution of shooting-rk4 combination when $y(s), y'(s)$ is not available.  

The results show that TIBA and TIBO are sensitive to volatile degree of underlying solution $y$ and the performance with $\theta=3\pi/2$ is better than that with $\theta=\pi/2$. Also, TIBA outperforms TIBO in general,  especially with $Dirichlet$ condition although TIBO is comparable to TIBA for other boundary conditions. Finally, TIBA and TIBO outperform $rk4$ significantly, especially the combination of shooting-rk4 becomes not dependable for $Mix_2$ where shooting method fails to identify required $y(s),y'(s)$. 
\begin{table}[htbp]
	\small
	\caption{The performance on non-homogeneous ODE}
	\begin{tabular}{lcccr}
		test type	& \multicolumn{1}{l}{$max(|y_{l}-y_{b}|)$} & \multicolumn{1}{l}{$max(|y_{o}-y_{b}|)$} & \multicolumn{1}{l}{$max|y_{rk4}-y_{b}|$} & $\theta$ \\ \hline\hline
		$Neumann$  & 4.7E-09 & 5.3E-09 & 8.6E-08 & $\pi/2$ \\
		$Neumann$  & 2.6E-08 & 2.6E-08 & 4.3E-06 & $3\pi/2$ \\
		$Dirichlet$     & 2.1E-12 & 7.0E-10 & 1.1E-08 & $\pi/2$ \\
		$Dirichlet$     & 3.1E-11 & 1.2E-09 & 9.1E-07 & $3\pi/2$ \\
		$Mix_1$   & 5.0E-10 & 2.4E-09 & 4.1E-08 & $\pi/2$ \\
		$Mix_1$   & 1.2E-08 & 1.2E-08 & 1.4E-06 & $3\pi/2$ \\
		$Mix_2$ & 2.3E-08 & 2.3E-08 & 1.0E-01 & $\pi/2$ \\
		$Mix_2$ & 1.7E-07 & 1.7E-07 & 1.0E-01 & $3\pi/2$ \\ \hline
	\end{tabular}%
	\label{tab:linear_non_homogeneous}%
\end{table}%

\newpage
\section{Summary}\label{sec:summary}
In this paper, we propose a  trigonometric-interpolation based approach (TIBA) algorithm to approximate solutions of second order nonlinear ODEs with general two-point linear boundary conditions.  It is comparable to another trigonometric-interpolation based approach developed in \cite{zou_tri_III} by leveraging analytic attractiveness of  trigonometric polynomial to capture the dynamics of $y,y',y''$  implied by ODE and boundary conditions. In a nutshell,   TIBA discrete globally ODE into a non-linear system.  TIBA is particularly attractive to copy with linear ODE where the solution can be solved directly without applying optimization process as in TIBO, which is not only more effective, but avoid dependence on initial values of target solution, which has significant impact on the TIBO's performance with a general mixed initial/boundary conditions.
In addition,  TIBA can copy with uniqueness and existence issues of solution properly. We have conducted numerical tests with four sets of boundary conditions and the results confirms that TIBA generally outperforms TIBO to handle linear ODE although TIBO should  be more effective to handle non-linear ODE. Finally,  TIBA's framework can be enhanced to solve second order linear Integro-differential equations with general initial/boundary mixed conditions as shown in \cite{zou_tri_fredholm} and \cite{zou_tri_volterra}. 
\appendix
\section{The derivation of Eq (\ref{eq:top}-\ref{eq:bottom})}\label{appendix:eq:UAV}
In this Appendix,  we derive Eq (\ref{eq:top}-\ref{eq:bottom}). Let $S, C$ be defined by Eq. (\ref{eq:SC}). Recall the following identity established in \cite{zou_tri}.
\begin{equation}\label{key_idenity}
\sum_{j=0}^{n-1} j \sin\frac{2\pi j k}{2n} =(-1)^{k+1}\frac n2 \cot\frac{\pi k}{2n}, \qquad 0<k<2n.
\end{equation}
We start with two special cases $u_0$ and $u_M$.  By Eq (\ref{eq:tildeu_d}), 
\begin{eqnarray}
	u_0 &=& a_0 I - \frac{b}{\pi } I' \times diag(1/K) B \nonumber\\
	&=& v_0 (\frac{2\pi }{bM^2}I'\cdot diag(K)\cdot SK -\frac{2\pi}{Mb} I' \cdot diag(K) S I -\frac{1}b) \nonumber\\
	&+&v_M(\frac{1}b   -\frac{2\pi }{bM^2}I'\cdot diag(K)SK  )   + \frac{2\pi }{bM}I' \cdot diag(K)SV \nonumber\\
	&=& v_0 (\frac{2\pi }{bM^2}\sum_{1\le i,j<M}ij\sin\frac{2\pi ij}{N}  -\frac{2\pi}{Mb} \sum_{1\le i,j<M}i\sin\frac{2\pi ij}{N} -\frac 1b)\nonumber\\
	&+&v_M(\frac 1b  -\frac{2\pi }{bM^2}\sum_{1\le i,j<M}ij\sin\frac{2\pi ij}{N})   + \frac{2\pi }{bM}   \sum_{1\le i,j<M}v_i j\sin\frac{2\pi ij}{N} \nonumber\\
	&=& v_0 (\frac{\pi }{bM}\sum_{1\le i <M}i (-1)^{i+1}\cot\frac{\pi i}N -\frac{\pi}{b} \sum_{1\le j<M} (-1)^{j+1}cot\frac{j\pi}{N} -\frac 1b)\nonumber\\
	&+&v_M(\frac 1b  -   \frac{\pi }{bM}\sum_{1\le i<M}i (-1)^{i+1}cot\frac{\pi i}{N} )   \nonumber\\
	&+& \frac{\pi }{b}   \sum_{1\le i<M}v_i (-1)^{i+1} \cot\frac{i\pi}{N}. \label{eq:u_0}
\end{eqnarray}
For the treatment of $u_M$, applying Eq \ref{key_idenity} to Eq (\ref{eq:tildeu_d}), we have
\begin{eqnarray}
	u_M &=& a_0 I - \frac{b}{\pi } I_a \cdot diag(1/K) B \nonumber\\
	&=& v_0 (\frac{2\pi }{bM^2}I_a\cdot  diag(K)SK -\frac{2\pi}{Mb} I_a \cdot diag(K) S I -\frac{1}b) \nonumber\\
	&+&v_M(\frac{1}b   -\frac{2\pi }{bM^2}I_a \cdot diag(K)SK  )   + \frac{2\pi }{bM}I_a \cdot diag(K)SV \nonumber\\
	&=& v_0 (\frac{2\pi }{bM^2}\sum_{1\le i,j<M}(-1)^iij\sin\frac{2\pi ij}{N}  \nonumber\\
	&-&\frac{2\pi}{Mb} \sum_{1\le i,j<M}(-1)^i i\sin\frac{2\pi ij}{N} -\frac 1b) \nonumber\\
	&+&v_M(\frac 1b  -\frac{2\pi }{bM^2}\sum_{1\le i,j<M}(-1)^i ij\sin\frac{2\pi ij}{N})   \nonumber\\
	&+& \frac{2\pi }{bM}   \sum_{1\le i,j<M} (-1)^jv_i j\sin\frac{2\pi ij}{N} \nonumber\\
	&=& v_0 (-\frac{\pi }{bM}\sum_{1\le i}i \cot\frac{\pi i}N +\frac{2\pi}{Mb} \sum_{1\le i,j<M}i\sin\frac{2\pi i(M-j)}{N} -\frac 1b) \nonumber\\
	&+&v_M(\frac 1b  +   \frac{\pi }{bM}\sum_{1\le i<M}i \cot\frac{\pi i}{N} )   - \frac{2\pi }{bM}   \sum_{1\le i,j<M} v_i j\sin\frac{2\pi (M-i)j}{N} \nonumber\\
	&=& v_0 (-\frac{\pi }{bM}\sum_{1\le i}i \cot\frac{\pi i}N + \frac{\pi}{b} \sum_{1\le j<M} (-1)^{M-j+1} \cot\frac{(M-j)\pi}{N}   -\frac 1b)\nonumber\\
	&+&v_M(\frac 1b  +   \frac{\pi }{bM}\sum_{1\le i<M}i \cot\frac{\pi i}{N})  \\
	&-& \frac{\pi }{b}   \sum_{1\le i<M} (-1)^{M-i+1}v_i \cot\frac{(M-i)\pi}{N}\nonumber\\
	&=& v_0 (-\frac{\pi }{bM}\sum_{1\le i}i \cot\frac{\pi i}N + \frac{\pi}{b} \sum_{1\le j<M} (-1)^{j+1} \tan\frac{j\pi}{N}   -\frac 1b)\nonumber\\
	&+&v_M(\frac 1b  +   \frac{\pi }{bM}\sum_{1\le i<M}i  \cot\frac{\pi i}{N}) - \frac{\pi }{b}   \sum_{1\le i<M} (-1)^{i+1}v_i \tan\frac{i\pi}{N}. \label{eq:u_M}
\end{eqnarray}
For other $u_i$ with $0<i<M$, we have by (\ref{eq:tildeu_d})
\begin{eqnarray}
	U(1:M-1) &=& a_0 I - \frac{b}{\pi } C \cdot diag(1/K) B \nonumber\\
	&=& a_0 I  - \frac{2a_1 \pi}{Mb} C \cdot diag(K) S I  \nonumber\\
	&-& \frac{2a_0\pi }{M^2}C\cdot diag(K)SK + \frac{2\pi }{bM}C\cdot diag(K)SV\nonumber\\
	&=& \frac{v_M-v_0}b I  - \frac{2v_0 \pi}{Mb} C \cdot diag(K) S I  \nonumber\\
	&-& \frac{2(v_M-v_0)\pi }{bM^2}C\cdot diag(K)SK \nonumber\\
	&+& \frac{2\pi }{bM}C\cdot diag(K)SV(1:M-1)\nonumber\\
	&=& v_0 (\frac{2\pi }{bM^2}C\cdot diag(K)SK -\frac{2\pi}{Mb} C diag(K) S I -\frac{1}bI) \nonumber\\
	&+&v_M(\frac{1}b I  -\frac{2\pi }{bM^2}C\cdot diag(K)SK  )   \nonumber\\
	&+& \frac{2\pi }{bM}C\cdot diag(K)SV(1:M-1). \label{eq:UMm1}
\end{eqnarray}
Let $c_i=C(i,:)$ and for $0<i<M$,  by (\ref{eq:UMm1}),
\begin{eqnarray}
	u_i &=& (v_M-v_0)(\frac1b - \frac{2\pi}{M^2b}c_i \cdot diag(K) \cdot S \cdot K) \nonumber\\
	&-& v_0\frac{2\pi}{bM}c_i \cdot diag(K) \cdot S \cdot I + \frac{2\pi}{bM}c_i \cdot diag(K) \cdot S \cdot V(1:M-1)\nonumber\\
	&=& v_0(\frac{\pi}{2b}\sum_{1\le k<M}(-1)^{i+k} \cot(k,i) \nonumber\\
	&-&  \frac{\pi}{2bM}\sum_{1\le k<M}(-1)^{i+k} k \cdot \cot(k,i) -\frac1b)   \nonumber\\
	&+& v_M(\frac{\pi}{2bM}\sum_{1\le k<M}(-1)^{i+k} k \cdot \cot(k,i) +\frac1b) \nonumber\\
	&-&\frac{\pi}{2b}\sum_{1\le k<M}(-1)^{i+k} \cot(k,i)v_k. \label{eq:ui}
\end{eqnarray}
One can see that Eq (\ref{eq:UAV}) represent the system \ref{eq:u_0}), (\ref{eq:u_M}), and (\ref{eq:ui}).


\begin{thebibliography}{99}
	\bibitem{ode_naga}
	{N. A. Gasilov},
	\textit{ On the existence and uniqueness of a solution to the boundary value problem for linear ordinary differential equations}, {Acta Math. Univ. Comenianae}, 93 (2024), pp {205-224}
	
	\bibitem{HBKeller}
	{Herbert B. Keller},
	\textit{Numerical Methods for Two-Point Boundary-Value Problems}, {Dover Publications, Inc. Mineola, New York},  2018
	
	\bibitem{ode_cuma_2}
	{S. Cuomo, A. Marasco},
	\textit{A numerical approach to nonlinear two-point boundary value problems for ODEs},
	{Computers and Mathematics with Applications}, 55 (2008),
	pp {2476-2489}
	
	\bibitem{losite}
	{José L. López, Ester Pérez Sinusía, Nico M. Temme},
	\textit{Multi-point Taylor approximations in one-dimensional linear boundary value problems},
	{Applied Mathematics and Computation}, {https://doi.org/10.1016/j.amc.2008.11.015},
	207 (2)(2009), pp {519-527}
	
	\bibitem{ode_17_chsa}
	{C. Chun, R. Sakthivel},
	\textit{Homotopy perturbation technique for solving two-point boundary value problems comparison with other methods}, {Comput. Phys. Commun}, 181 (2010), pp {1021–1024}
	
	
	\bibitem{ode_Adomian}
	{G. Adomian},
	\textit{ A review of the decomposition method in applied mathematics}, {J. Math. Analysis and Applications}, 135 (1988), pp {501-544}
	
	
	\bibitem{ode_19_gecu}
	{F. Geng, M. Cui},
	\textit{A novel method for nonlinear two-point boundary value problems: combination of ADM and RKM},
	{Appl. Math. Comput}, 217 (2011), pp {4676–4681}
	
	\bibitem{ode_24_kare}
	{A.S.V. Ravi Kanth, Y.N. Reddy},
	\textit{Cubic spline for a class of singular two-point boundary value problems}, {Appl. Math. Comput.}, 170 (2005), pp {733–740}
	
	
	\bibitem{ode_zhlish}
	{Y. Zheng, Y. Lin, Y. Shen},
	\textit{A new multiscale algorithm for solving second order boundary value problems}, {Applied Numerical Mathematics}, 156 (2020), pp { 528–541},
	
	\bibitem{ode_3}
	\textit{Topological methods for some boundary value problems}, {H.B. Thompson},
	{J. Comput. Math. Appl.}, 42 (2001), pp {487–495},
	
	\bibitem{ode_8}
	{M. Al-Smadi, O. Abu Arqub, N. Shawagfeh, et al.},	
	\textit{Numerical investigations for systems of second-order periodic boundary value problems using reproducing kernel method}, {Appl. Math. Comput.},291 (2016),pp {137–148}
	
	
	\bibitem{ode_22}
	{W. Jiang, M. Cui},
	\textit{ Solving nonlinear singular pseudoparabolic equations with nonlocal mixed conditions in the reproducing kernel space}, {Int. J. Comput. Math.}, 87 (2010), pp {3430–3442},
	
	\bibitem{ode_36}
	{Z.H. Zhao, Y.Z. Lin, J. Niu},
	\textit{Convergence order of the reproducing kernel method for solving boundary value problems}, {Math. Model. Anal.}, 21 (2016), pp {466–477}
	
	\bibitem{ode_bumozh13}
	{S. Bushnaq, S. Momani, Y. Zhou},
	\textit{A reproducing kernel Hilbert space method for solving integro-differential equations of fractional order},
	{J. Optim. Theory Appl.}, 2014 (2013), pp {96–105}
	
	\bibitem{ode_agwoli}
	{R.P. Agarwal, F.H. Wong, W.C. Lian},
	\textit{Positive Solutions for Nonlinear Sigular Boundary Value Problems}, {Applied Math Letters}, 12 (1999), pp {115-120},
	
	\bibitem{ode_wo1}
	{F.H. Wong},
	\textit{Existence of positive solutions of sigular boundary value Problems}, {Nonlinear Analysis}, 21 (1993), pp {397-406},
	
	\bibitem{ryts}
	{Victor S. Ryaben'kii, Semyon V. Tsynkov},
	\textit{A Theoretical Introduction to Numerical Analysis}, {Chapman and Hall/CRC}, 2007
	
	\bibitem{zou_tri}
	{X. Zou}, 
	\textit{On Trigonometric Approximation and Its Applications  \qquad\qquad\qquad\qquad\qquad}
	{https://arxiv.org/pdf/2505.02330.}
	

	\bibitem{zou_tri_III}
	{X. Zou}, \textit{Trigonometric Interpolation Based Optimization for Solving Second Order Non-Linear ODE  with Mixed Boundary Conditions}, {https://arxiv.org/pdf/2504.19280}
	
	\bibitem{zou_tri_V}
	{X. Zou}, 
	\textit{Trigonometric Interpolation Based Optimization for High Order Non-Linear ODE}, To be appear.
	
	\bibitem{zou_tri_fredholm}
	{X. Zou}, 
	\textit{Trigonometric Interpolation Based Approach for Second Order Fredholm Integro-differential equations}, To be appear.
	
	\bibitem{zou_tri_volterra}
	{X. Zou}, 
	\textit{Trigonometric Interpolation Based Approach for Second Order Volterra Integro-differential equations}, To be appear.

\end{thebibliography}
\end{document}